\documentclass [12pt,twoside,a4paper] {article}
\usepackage{graphicx}
\usepackage{latexsym}
\usepackage{amssymb}
 \newfont{\bit}{cmbxti10 scaled \magstep1}  

  \newcommand {\kj} {\v{c}}
  \newcommand {\tj} {\'{c}}

  \newcommand {\ra} {\rightarrow}
  \newcommand {\lra} {\leftrightarrow}
  \newcommand {\longra} {\longrightarrow}

\author {Boris \v{C}ulina}
\title  {The  Concept of Truth }
\date   {Zagreb,\ 1999.}

\begin {document}
  \bibliographystyle {alpha}  
%
%
\def\rmdj {d\llap{\raise 1.22ex\hbox
  {\vrule height 0.09ex width 0.315em}\kern 0.04em}}
\def\sldj {d\llap{\raise 1.22ex\hbox
  {\vrule height 0.09ex width 0.265em}}\rlap{\raise 1.22ex\hbox
  {\vrule height 0.09ex width 0.05em}}}
\def\itdj {d\llap{\raise 1.22ex\hbox
  {\vrule height 0.09ex width 0.2em}}\rlap{\raise 1.22ex\hbox
  {\vrule height 0.09ex width 0.06em}}}
\def\bfdj {d\llap{\raise 1.16ex\hbox
  {\vrule height 0.126ex width 0.308em}\kern 0.04em}}
\def\ttdj {\rlap{\kern 0.17em\raise 1.1ex\hbox
  {\vrule height 0.09ex width 0.295em}}d}
\def\scdj {\rlap{\kern 0.04em\raise 0.57ex\hbox
  {\vrule height 0.09ex width 0.20em}}d}
\def\sfdj {d\llap{\raise 1.22ex\hbox
  {\vrule height 0.10ex width 0.3em}\kern 0.02em}}

\def\dj{\ifcase\fam \rmdj \or \or \or
  \or \itdj \or \sldj \or \bfdj \or \ttdj \or \sfdj \or \scdj \else \rmdj \fi}

\def\rmDj {\rlap{\kern 0.05em\raise 0.76ex\hbox
  {\vrule height 0.10ex width 0.28em}}D}
\def\slDj {\rlap{\kern 0.1em\raise 0.76ex\hbox
  {\vrule height 0.1ex width 0.28em}}D}
\def\itDj {\rlap{\kern 0.145em\raise 0.76ex\hbox
  {\vrule height 0.1ex width 0.274em}}D}
\def\bfDj {\rlap{\kern 0.044em\raise 0.72ex\hbox
  {\vrule height 0.126ex width 0.287em}}D}
\def\ttDj {\rlap{\kern 0.02em\raise 0.67ex\hbox
  {\vrule height 0.105ex width 0.20em}}D}
\def\scDj {\rlap{\kern 0.08em\raise 0.73ex\hbox
  {\vrule height 0.12ex width 0.24em}}D}
\def\sfDj {\rlap{\kern 0.02em\raise 0.727ex\hbox
  {\vrule height 0.126ex width 0.26em}}D}

\def\Dj{\ifcase\fam \rmDj \or \or \or
  \or \itDj \or \slDj \or \bfDj \or \ttDj \or \sfDj \or \scDj \else \rmDj \fi}

\def\int {\intop\limits}
\def\Intnorm#1#2
  {\int_{\lower1.5ex\hbox{$#1$}}^{
  \raise.5ex\hbox{$#2$}}}
\def\Intsmall#1#2
  {\int_{\hbox{\small\lower1.5ex\hbox{\small$#1$}}}^{
  \hbox{\small\raise.5ex\hbox{\small$#2$}}}}
\def\Intfoot#1#2
  {\int_{\hbox{\footnotesize\lower1.5ex\hbox{\footnotesize$#1$}}}^{
  \hbox{\footnotesize\raise.5ex\hbox{\footnotesize$#2$}}}}

   \maketitle
   \voffset=0in
 
   \pagenumbering{arabic}

\begin{flushright}``{\em The ghost of the Tarski hierarchy is still with us.}''\\
 Saul Kripke \end{flushright}
 
   {\bf Abstract} \ On the basis of elementary thinking about language functioning, 
a solution of truth paradoxes is given and a corresponding
semantics of a truth predicate  is founded. It is shown that it is precisely the two - valued
description of the maximal intrinsic fixed point of the strong Kleene three - valued
semantics. 

\section{Analysis of  the truth concept and the informal descripton of 
the solution}\mark{}{}

  Roughly, by the ``{\bf classical language }'' will be meant every language which is modelled
upon the everyday language of declarative sentences. An example is the  
standard mathematical language which is basically an everyday language supported by the 
symbolisation process and by the  mechanism of variables.
Due to definitness, the language of the first order logic,  which has an explicit and precise 
description of form and meaning, will be considered. By the ``language'' will be meant
an interpreted language, a language form together with an interpretation. 

   Besides a  formal (grammatical)  structure and  an internal meaning structure,  a language
has an external meaning structure too, a connection between language forms and external objects
which constitute the subject of the language. 
 For the classical language there are assumptions that there are objects which the language mentions,
that every name is a name of some object, that to  every functional and relational simbol 
  an operation  or a relation between objects is associated, and that every atomic sentence is
true or false, depending on if ``it is'' or `` it is not''   its content.
 These assumptions have grown from everyday use of language where we are accustomed to their
fulfilment, but there are situations when they are not fulfiled.
The Liar paradox and other  paradoxes of truth are witnesses of  such situations. 
They are  the results of a  tension between implicitly accepted assumptions on the language 
and their unfulfilment.

  Let's investigate the sentence $L$ (the {\bf Liar}): 

  $L$: $L$ is a false sentence. (or ``This sentence is false.'')  
  
  using the  everyday  understanding of truth and language , to investigate  truth of  
 $L$ we must investigate  what it says.
But it says precisely about its own truth, and  in a contradictory way. If we 
asume it is true, then it is true what it says --- that it is false. But if we assume
it is false, then it is false what it says, that it is false, so it is true. Therefore,
it is a selfcontradictory sentence. But what is even more important  is a paradoxical feeling
that we can't determine its truth value. The same paradoxality, but whithout contradiction, 
emerges during the  investigation of the following sentence $I$ (the {\bf Truthteller}):

$I$: $I$ is a true sentence. (or ``This sentence is true.'')

  Contrary to the Liar to which we can't associate any truth value, to this sentence we can 
 associate   the truth as well as the falshood with equal mistrust.
There are no additional specifications which would make a choice between the two possibilities,
not because we haven't enough knowledge but principally. Therefore, we can't associate a truth
value to these sentence, neither. 

  There are various analysis and solutions which will not be considered here
(good surveys   can be found in \cite{mar},\cite{hel},\cite{vis},\cite{she}).  
This analysis begins with a basic intuition that the previous sentences are meaningful
(because we understand well what they say, even more, we used that in the unsuccessful determination
of their truth values), but they witness the failure of the classical procedure for the truth value 
determination in some ``extreme'' situations. Paradoxality emerges from a confrontation
of the implicit assumption of the success of the procedure and the discovery of the failure.   
A basic assumption about the classical language is that every sentence is true ($\top$)
or false ($\bot $). Truth values of more complex sentences are determined according to truth values of
simpler components in a way determined by the internal 
semantics of the language. 
To visualize better this semantical relationship of sentences we will imagine them as nodes of a
graph $S$ and we will draw arrows from a sentence to all sentences on which its truth value 
depends. Because of definitness we will consider sentences of an interpreted first order language.
Due to simplicity we will assume that for every object $a$ of a domain $dom(L)$ of the 
language $L$ there is a closed term $\bar{a}$ which names it. The  logical vocabulary of the language
is standard and it consists of connectives $\neg ,\wedge ,\vee ,\ra ,\lra$ and quantors 
$\forall ,\exists$. Arrows of the semantical graph $S$ are defined by the recursion on inductive structure
of sentences. Instead of  strict definitions ``pictures'' of tipical nodes of the graph will be 
shown:


\begin{center}
	\includegraphics[width=\linewidth]{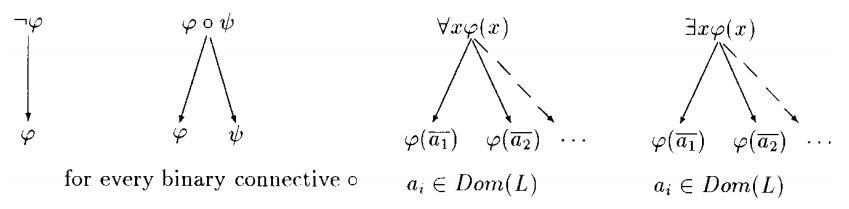}
\end{center}

 The truth values of the sentences will be described as a function $I: S\longrightarrow \{\top ,\bot\}$, 
which associates the truth value, $\top$ or $\bot$, to every sentence.
   
  The internal semantics of the language describes  determination of the truth value of a 
compound sentence in terms of the truth values of simpler sentences to which it shows (by  arrows 
of the semantical graph). The description consists of standard conditions on the truth function $I$:

\begin{enumerate}
\item $\displaystyle I(\neg\varphi )=\left\{\begin{array}{ll}
                                             \top & for\ I(\varphi )=\bot \\
                                              \bot & for\ I(\varphi )=\top 
                                             \end{array}\right. $
\item $\displaystyle I(\varphi \wedge \psi )=\left\{\begin{array}{ll}
                                             \top & for\ I(\varphi )=\top\ and\  I(\psi )=\top\
                                               (both\ are\ true)\\
                                              \bot &  for\ I(\varphi )=\bot\ or\  I(\psi )=\bot\
                                               (at\ least\ one\ is\ false) 
                                             \end{array}\right. $

\item $\displaystyle I(\varphi \vee \psi )=\left\{\begin{array}{ll}
                                             \top & for\ I(\varphi )=\top\ or\  I(\psi )=\top\
                                               (at\ least\ one\ is\ true)\\
                                              \bot &  for\ I(\varphi )=\bot\ and\  I(\psi )=\bot\
                                               (both\ are\ false) 
                                             \end{array}\right. $

\item $\displaystyle I(\varphi \ra \psi )=\left\{\begin{array}{ll}
                                             \top & for\ I(\varphi )=\bot\ or\  I(\psi )=\top\
                                                   \\
                                              \bot &  for\ I(\varphi )=\top\ and\  I(\psi )=\bot\
                                               
                                             \end{array}\right. $

\item $\displaystyle I(\varphi \lra \psi )=\left\{\begin{array}{ll}
                                             \top & for\ I(\varphi )=I(\psi )\
                                               (both\ are\ true\ or\ both\ are\ false)\\
                                              \bot &  for\ I(\varphi )\neq I(\psi )\
                                               (one\ is\ true\ and\ another\ is\ false) 
                                             \end{array}\right. $

\item $\displaystyle I(\forall x\varphi (x))=\left\{\begin{array}{ll}
                                             \top & if\ \forall a\in dom(L)\ I(\varphi (\bar{a}))=\top \\
                                              \bot & if\ \exists a\in dom(L)\ I(\varphi (\bar{a}))=\bot 
                                             \end{array}\right. $

\item $\displaystyle I(\exists x\varphi (x))=\left\{\begin{array}{ll}
                                             \top & if\ \exists a\in dom(L)\ I(\varphi (\bar{a}))=\top \\
                                              \bot & if\ \forall a\in dom(L)\ I(\varphi (\bar{a}))=\bot
                                             \end{array}\right. $

\end{enumerate}  

  According to  these conditions, to determine the truth value of a given sentence we must 
investigate the truth values of all sentences  which it shows, then eventually, for the same reasons,
the truth values of the sentences  which these sentences show, and so on. Every such path along the arrows
of the graph leads to atomic sentences (because the complexity of  sentences decreases along 
the path) and the  truth value of the initial sentence is completely determined by  the truth 
values of atomic sentences which it hereditary shows. In common situations language doesn't 
talk about the truth values of its own sentences, so  the truth values of its atomic sentences don't
depend on the truth values of some other sentences. They are leafs of the semantic graph --- there
are no arrows from them leading to other sentences. To investigate their truth values we must investigate
external reality they are talking about. For example, for the simplest atomic sentence $P(\bar{a})$
to investigate its truth value we must see whether the object $a$ has the property $P$. The assumption 
of the classical language is that it is or it is not the case --- $P(\bar{a})$ is true or false.
It is fulfiled in standard situations, whether effectively or principally. Therefore, every atomic
sentence has a definite truth value, so the procedure of determination of the truth value of
 every sentence also
gives a definite truth value, $\top$ ili $\bot$. Formally, it is secured by the recursion
principle which says that there is a unique function $I: S\longrightarrow \{\top ,\bot\}$ with 
 values on atomic sentences being identical to externally given truth values, and it obeys
previously displayed classical semantic conditions.  
 
 But the classical situation can be (and it is) destroyed when atomic sentences talk about the 
truth values of other sentences. Then there are arrows from atomic sentences to other sentences
along which we must continue to investigate the truth value of the initial sentence. The simplest such
a situation is when language contains the {\bf truth predicate} $T$ by means of which it can talk about
the truth values of its own sentences. Then language has atomic sentences of the form $T(\bar{\varphi})$
with the meaning ``$\varphi$ is a true sentence''. The truth conditions for $T(\bar{\varphi})$ 
are  part of the {\em internal semantics of the language}, as there are for example the truth 
conditions on $\varphi\wedge\psi$. They don't depend on the external world but on the truth value of
 the sentence $\varphi$ by a  {\em logical sense} we associate to the truth predicate $T$---
 we consider $T(\bar{\varphi})$ to be true when $\varphi$ is true, and to be false when
$\varphi$ is false.
  
 So, in the case of  presence of the truth predicate $T$ there are new arrows in the graph


\begin{center}
	\includegraphics[height=3cm]{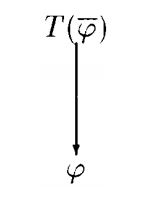}
\end{center}

 and a new condition on the truth function:

\

 $\displaystyle I(T(\bar{\varphi}))=\left\{\begin{array}{ll}
                                             \top & for\ I(\varphi )=\top \\
                                              \bot & for\ I(\varphi )=\bot 
                                             \end{array}\right. $
  
\
 
Now, to investigate the truth value  of a sentence  it is not sufficient to reduce the problem
to atomic sentences in general, but we must again continue the ``voyage'' upon arrows  to more
complex sentences. Because of  the possible ``circulations'',  there is nothing  to insure 
the success of the procedure. Truth paradoxes just witness such situations. Three illustrative
examples follow


\begin{center}
	\includegraphics[height=3cm]{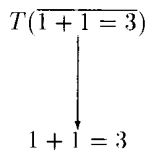}
\end{center}

 The procedure of the truth value determination has stopped on the atomic sentence for which we know 
 is false, so $T(\overline{1+1=3})$ is false, too.

The Liar: For $L:T(\overline{\neg L})$ we have


\begin{center}
	\includegraphics[height=3cm]{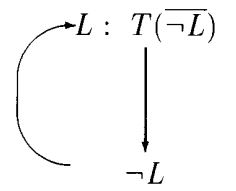}
\end{center}

  But now the procedure of  the truth value determination has failed because the conditions for
the truth function can't be fulfiled. Truth value of  $T(\overline{\neg L})$ depends on
truth value of $\neg L$ and this again on $L:T(\overline{\neg L})$ in a way  which 
 is impossible to obey.

The Truthteller: For $I:T(\bar{I})$ we have


\begin{center}
	\includegraphics[height=2cm]{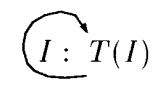}
\end{center}

 Now, there are, as we have already seen, two possible assignings of truth values to the sentence 
$I$. But this multiple fulfilment we must consider as a failure of the classical procedure, too,
because it assumes to establish a unique truth value for every sentence.

Paradoxes  emerge just because the classical procedure of the truth value determination sometimes 
doesn't give a classically assumed (and expected) answer. As previous examples show such assumption
is an unjustified generalization from common situations to all situations. We can preserve the 
classical procedure, also the internal semantic structure of the language. But, we must reject 
universality of the assumption of its success. The awareness of that transforms paradoxes
to normal  situations inherent to the classical procedure. I believe this is the solution of paradoxes.
But, there remains the solving of another significant question  --- how to insure a success of the truth
value determination procedure which is crucial for the validity of the classical logic, and in the
same time to preserve the internal semantic structure of the language. Certainly,  prohibition of a 
language which talks of its own truth can't be considered as a satisfactory solution, nor
the hierarchy of languages in which every language can talk only about truth of sentences which
belong to the language below it in the hierarchy. Although circularity is a substantial part of 
paradoxical situations its rejection is a too rough solution which impoverishes language 
inacceptabily. As Kripke showed in \cite{kri} circularity is deeply present in an everyday 
language use not only in an unavoidable way but also in a harmless way, and only in some 
extreme situations it leads to paradoxes. Kripke showed it on examples which involve
external meaning structure of language (``empirical facts''), but the same occurs in 
internal meaning structure, too. Neither there circularity leads necessary to paradoxes,
as the following example shows.

 Let's determine truth value of the sentence the {\bf Logician}:

$Log:\ T(\overline{Log})\ \vee\ T(\overline{\neg Log})$ (This sentence is true or false)

       Semantical dependencies are the following:


\begin{center}
	\includegraphics[height=5cm]{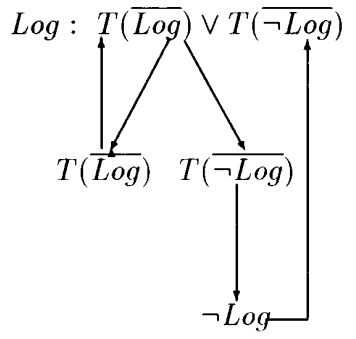}
\end{center}

If $Log$ were false then, by the truth conditions, $T(\overline{\neg Log})$ would be false,
$\neg Log$ would be false too, and finally $ Log$ would be true. Therefore, such valuation of
the graph is impossible. But if we assume that $ Log$ is true, the truth conditions  generate
a unique consistent valuation. Therefore, the truth determination procedure gives the unique answer
--- that  $ Log$ is true.

Kripke showed in \cite{kri} that circumstances which lead to paradoxes cannot be isolated on a 
sintactical level, but an intervention in the semantic language structure is  neccessary.
The intervention is here made  in the following way. The primary classical semantics is preserved,
so the classical procedure of truth value determination is preserved too, but the wrong classicall
 assumption of its total success is rejected. The rejection doesn't change the meaning of the classical
conditions on the truth function, because they are stated in a way independent of the assumption
that the function is everywhere defined. Their functioning in the new situation is illustrated in
the following sentence:

$L\ \vee\ \overline{0=0}$

 On the classical condition for the connective $\vee$ this sentence is true precisely when at least one
of the basic sentences is true. Because  $\overline{0=0}$ is true consequently the total
sentence is true regardless of the fact that $L$ hasn't the truth value. Equally, if we apply the
truth value condition on the connective  $\wedge$ to the sentence

$L\ \wedge\ \overline{0=0}$

  the truth value will not be determined. Namely, for the sentence to be true  both basic
sentences must be true, and it is not fulfiled. For it to be false at least one basic sentence must
be false and this also is not fulfiled. So, nonexistence of the truth value for $L$ leads to nonexistence
of the truth value for the whole sentence.

 Classical truth value conditions specify the truth value of a compound sentence in terms of
truth values of its direct components regardless whether they have truth values or not. The lack
of some truth value  may lead, but does not have to, to the lack of the truth value of the compound sentence.
It is completely determined by the classical meaning of the construction of a sentence and by
the basic assumption that all sentences are considered meaningful regardless of the truth value. 

 Therefore, some sentences, although meaningful,  valued by the classical conditions have not 
the truth value, because the conditions do not give them a unique truth value. This leads to the
{\bf partial two - valued semantics of the language}. Where the procedure gives a unique truth value, 
truth or falshood, we accept it, where it fails because it does not give any truth value or permits 
both values, the sentence remains without the truth value. This kind of semantics can be described as the
{\bf three - valued semantics of the language} --- simply the failure of the procedure will be 
declared as a third value $\mid$ ({\bf undetermined}). It has not any additional philosophical charge.
It is only a convenient technical tool for the description.  

  But this semantics is not accepted here as the final semantics of the language. A decidable 
reason for  the rejection is the opinion  that  the two - valued semantics
is  natural  to human kind and that every other semantics can be reduced to the two - valued  by an
appropriate modelling (a confirmation for the thesis is that descriptions of all semantics are 
two - valued). To remain on three - valued semantics would mean that the logic would not be classical, 
the one we are accustomed to. Concerning the truth predicate itself, it would imply the  
preservation of its classical 
logical sense in  the two - valued part of the language extended by the ``silence'' in the
part where the classical procedure fails. Although in a metadescription $T(\bar{\varphi})$ has 
the same  truth value (in the three - valued
semantic frame) as $\varphi$, that semantics is not the more  initial one (although it extends it)
nor it can be expressed in the language itself (because the language is silent about the third value,
or better said, the third value is the reflection in the metalanguage of the silence in the language).
So the expressive power of the language is weak. For example, the Liar is undetermined. Although
 we have easily
said it  in metalanguage we cannot express in the language $L$ itself, because, as it has already been
said (in metalanguage), the Liar is undetermined. Not only that the ``zone of silence'' is
unsatisfactory because of the previously stated reasons (it leads to the three - valued logic, it loses the 
primary sense of the truth predicate and it weakens the expressive power of the language), but it can 
be interrupted by a natural {\em additional valuation} of the sentences which emerges from 
recognising the failure of the classical procedure. This point  will be illustrated on the example of 
the Liar. On the
intuitive level of thinking, by recognising the Liar is not true nor false we state that it is 
undetermined. So, it is not true what it claims --- that it is false. Therefore, the Liar is false.
But this does not lead to  restoring of the contradiction because  a {\em semantical
shift} has happened from the primary partial two - valued semantics (or three - valued semantics) toward its two - valued 
description,  which merely extends it in the part where it is not determined. Namely,
the Liar talks of its own truth in the frame of the primary semantics, while the last valuation is      
in the frame of the final semantics. The falsehood of the Liar in the final semantics doesn't
mean that it is true what it says (that it is false) because the semantical frame is not the same. It means that
it is false what it talks of its ownprimary semantics (that it is false in the primary semantics). It follows that it is not false in
its primary semantics. But, it cannot be true in the primary semantics because then it would be
true in the final semantics (which only extends the primary  where it fails). Therefore
it is undetermined in the primary semantics. So, not only have we gained a contradiction, but we
 also have received
another information about the Liar.

 It is easy to legalize this intuition. Using the truth predicate, the language talks about its primary
semantics. The classical procedure and the classical meaning of the truth predicate determines
its primary semantics, which is, due to the failures of the procedure,  a partial
twovalued semantics (= a threvalued semantics). {\em But the description of the primary semantics itself
is its natural extension to the final two - valued semantics.} Therefore, the final semantics of the
language has for its subject precisely the primary semantics of the language which it extends furthermore 
 in the part where it is silent  using the informations about the silence. The transition can be 
described easily on the semantical graph. To get the final valuation from the primary valuation
we must  revaluate only  the atomic sentences of the form $T(\bar{\varphi})$. Such  sentences have the 
{\em same meaning in both semantics} --- that $\varphi$ is true in the primary semantics, but the
{\em truth conditions are not the same}. While in the primary semantics the truth conditions
for $T(\bar{\varphi})$ are classical  (the truth of $T(\bar{\varphi})$ means the
truth of  $\varphi$, the falsehood of $T(\bar{\varphi})$ means the falsehood of $\varphi$), 
in the final semantics it is not so. In it the truth of $T(\bar{\varphi})$ means that $\varphi$
is true in the primary semantics, and  falsehood of $T(\bar{\varphi})$ means that $\varphi$
is not true in the primary semantics. It does not mean that it is false in the primary semantics,
but that it is false or undetermined. So,  formally looking,  in the final semantics 
$T(\bar{\varphi})$ inherits  truth from the primary semantics, while other values transform
to  falsehood.
 
 We can see  best that this is a right and a complete description of the valuation in the
 primary semantics  by introducing predicates for other truth values in the primary valuation:

 $F(\bar{\varphi})$ (= $\varphi$ is false in the primary semantics) $\lra$ $T(\overline{\neg\varphi})$ 
   
 $U(\bar{\varphi})$ (= $\varphi$ is undetermined in the primary semantics) $\lra$
 $\neg T(\bar{\varphi})$ i $\neg F(\bar{\varphi})$

 According to the truth value of the sentence $\varphi$ in the primary semantics we determine 
which of the previous sentences are true and which are false. For example, if  $\varphi$ is false in
the primary semantics then $F(\bar{\varphi})$ is true while others
($T(\bar{\varphi})$ i $U(\bar{\varphi})$) are false.

 Once  the final two - valued valuations of atomic sentences are determined in this way,
  valuation of every sentence is determined  by means of the classical conditions and  the 
principle of recursion. This valuation not only preserves the primary logical meaning of the truth predicate
(as the truth predicate of the primary semantics) but it also coincides with the primary valuation
 where it is determined. Namely, if $T(\bar{\varphi})$ is true in the primary semantics then 
$\varphi$ is true in the primary semantics, so $T(\bar{\varphi})$ is true in the final semantics.
If $T(\bar{\varphi})$ is false in the primary semantics then $\varphi$ is false in the primary 
semantics, so $T(\bar{\varphi})$ is false in the final semantics. Since the truth conditions for
compound sentences are the  same in both semantics this coincidence spreads through all sentences which
have determined value in the primary valuation. Therefore $T(\bar{\varphi})\ra \varphi$ and
 $F(\bar{\varphi})\ra \neg\varphi$ are true sentences in the final semantics.

 Having in mind this kind of double semantics of the language, we can easily solve all  truth paradoxes. 
On an intuitive level we have already done it for the Liar. To distinguish inside which semantic frame we use 
a certain term we will put prefix ``p'' for the primary semantics and prefix ``f'' for the final
semantics. In that way we will distinguish for example ``f-falsehood'' and ``p-falsehood. The form
of the solution is always the same. A paradox in the classical thinking means that the
truth value of a sentence is undetermined in the primary semantics. But, then it becomes an information in the final semantics
with which we can conclude the truth value of the sentence in the final semantics.

First, let's investigate  the situations which lead to the contradiction like the Liar. Of such kind is,
for example, the {\bf Strong Liar} $LL:\neg T(\overline{LL})$ (``This sentence is not true'').
In the naive semantics it leads to a contradiction in the same way as the Liar, because there
``not to be true'' is the same as ``to be false''.  Recognising a failure of the classical 
procedure, we continue to think in the final semantics and state that it is p-undetermined. So, 
it is not p-true. But, it claims just that, so it is f-true. Therefore, we conclude that the Strong
Liar is undetermined in the primary semantics and true in the final semantics. It is interesting  that
 the whole argumentation can be done directly in the final semantics,  not indirectly by stating
the failure of the classical procedure. The argumentation is the following. If $LL$ were f-false, then
it would be f-false what it said --- that it is not p-true. So, it would be p-true. But, it means
(because the final semantics extends the primary one) that it would be f-true and it is a contradiction 
with the assumption. So, it is f-true. This statement does not lead to a contradiction but to an
additional information. Namely, it follows that what it talks about is f-true --- that it is not 
p-true. So, it is p-false or p-undetermined. If it were p-false it would be f-false too, and this is
a contradiction. So, it is p-undetermined. Therefore, although the Liar and the Strong Liar
are both p-undetermined, the latter is f-true while the former is f-false.

 Let's analyse in the same way {\bf Curry's paradox} $C:T(\bar{C})\ra l$ (``If this sentence is true
then $l$''), where $l$ is any false statement. On the intuitive level if $C$ were false then 
$T(\bar{C})$ is true, and so is $C$ itself, and it is a contradiction. If $C$ was true then the
whole conditional ($C$) and its antecedent $T(\bar{C})$ would be true, and so the consequent
$l$ would be true, which is impossible with the choice of $l$ as a false sentence. Therefore we conclude
in the final semantics that $C$ is p-undetermined, and so it is f-true (because the antecedent is
f-false). The argumentation can  also be completely translated in the final semantics as follows.  
Namely, if $C$ were f-false then the antecedent would be f-true. It means that $C$ would be p-true 
and therefore f-true (by the accordance of two semantics), and it is a contradiction. So, $C$ is
f-true.  From it we conclude  that $T(\bar{C})$ is f-false or  $l$ is f-true. Because $l$ is f-false
it follows that $T(\bar{C})$ is f-false, so $C$ is not  p-true. It is therefore p-false or
p-undetermined. If it were p-false it would be f-false, so it is p-undetermined.

 In the same way other truth paradoxes, which lead to contradiction  on the intuitive level, lead
to positive argumentation in the final semantics. But, the situation is  different with paradoxes
which do not lead to contradiction, which permit more valuations, like the {\bf Truthteller}. Its 
analyses gives that it is p - undetermined. It implies that it is not p - true which means that  
($I:T(\bar{I})$)  it is not $I$.
So, $I$ is f-false. But, this thinking cannot be translated directly into the final semantics.
 The argumentation formulated in the final semantics do not give the answer as well as in the primary 
semantics. It is necessary to investigate primary valuations of the semantical graph. Of course,
if we enrich the language with the description of  semantical graphs and truth valuations then
it is possible to translate the intuitive argumentation.

 Through this kind of modeling, the truth predicate of the primary semantics is described by the final
semantics of the language. Of course, it does not coincide with the truth predicate of 
the final semantics. But the goal was not to describe the predicate. Moreover, it emerges from the
description of the primary truth predicate. Being at  the same time an extension of the
primary predicate, it describes itself partially, but not completely. In that sense the ghost
of the Tarski hierarchy is still with us. For some it is an evil ghost because it does not
permit the complete description of the truth predicate of the final semantics. For the author
it is a good ghost because thanks to him the truth predicate of the primay semantics is 
completely described.

\section{Formal description of the solution}\mark{}{}

  Let $L$ be an interpreted first order language with a domain $D$. It will be
permited that $D$ can be an empty set. Then $L$ reduces to its logical vocabulary. Also, because of
simplicity, we will assume that for every object $a\in D$ of the language there is a closed term
$\bar{a}$ which names it.

  We will extend $L$ to the language $LT$ which will talk additionally  about the truth of its own
sentences. Along with objects of the language $L$ its domain will contain its own sentences, too.
Its vocabulary will contain the predicate $S$ (= ``to be a sentence'') which will distinguish
 sentences from  other objects and the predicate $T$ (= ``to be a true sentence'') which will
 describe the truth of the sentences. Every sentence $\varphi$ will have its own name $\overline{\varphi}$,
but there will be also special names for sentences. Giving to them  suitable denotations
we will achieve intended selfreferences. For example, the constant $\overline{I}$ will be interpreted
as a name of the sentence $T(\overline{I})$, and so we will construct the {\bf Truthteller} in the
language.

{\bf Vocabulary} of the language $LT$ consists of the vocabulary of the language $L$ together with 
new simbols --- unary predicates $S$ and $T$,  {\bf sentence constants} $\overline{I}$, $\overline{L}$, $\overline{\neg L}$,
$\overline{LL}$, $\ldots$ and the special operator  $\bar{\ }$.

  A set of {\bf terms} $TLT$ and a set of  {\bf formulas} $FLT$ of the language $LT$ are defined
 as the smallest sets which
satisfy all conditions for terms and formulae of the language $L$ and additional conditions:

\begin{enumerate}
\item Sentence constants $\overline{I}$, $\overline{L}$, $\overline{\neg L}$,
      $\overline{LL}$, $\ldots$ are terms
\item If $t$ is a term then $S(t)$ and $T(t)$ are formulae
\item If $\varphi$  is a formula then $\overline{\varphi}$ is a term. 
\end{enumerate}

 A set of {\bf free variables} of  a term or a  formula are defined by standard recursive 
conditions plus one more condition --- that free variables of the term  $\overline{\varphi}$ 
are precisely free variables of a formula $\varphi$. {\bf Sentences} of the language $LT$ are
closed formulae of the language. A set of them will be marked $SLT$.

 The {\bf interpretation}({\bf model}) of the language $LT$ is given in the following way. 
Domain $DLT$ consists of all objects from domain $D$ together with all sentences from $LT$:
$DLT=D\cup SLT$. All predicates of the language  $L$ are extended inside new domain in  such manner
that they give falsehood if at least one argument is outside domain $D$, and functions of the language $L$
are extended in  such manner that they give some constant value, let's say the sentence  $T(\overline{I})$,
if at least one argument is out of the domain $D$. 

  New simbols are interpreted in the following way. The simbol $S$ has the meaning ``to be a sentence'',
that is, it is interpreted by the set of sentences of $LT$. The simbol $T$ will be the truth 
predicate of the primary semantics, in other words it will have the meaning ``to be a true sentence
in the primary semantics''  once we define what  the primary semantics is. In the primary semantics it will be 
achieved in such manner that $T$ will be introduced as a logical symbol (like $\wedge$ for 
example) with the classical truth conditions, and in the final semantics it will be achieved directly
interpreting it as the set of true sentences in the primary semantics. The sentence constants will be 
interpreted as names of appropriate sentences:

\begin{enumerate}
\item $\overline{I}$ is the name of a sentence $T(\overline{I})$
\item $\overline{L}$ is the name of a sentence $T(\overline{\neg L})$
\item $\overline{\neg L}$ is the name of a sentence $\neg T(\overline{\neg L})$ 
\item $\overline{LL}$ is the name of a sentence $\neg T(\overline{LL})$
\end{enumerate}
\hspace{1cm} and so on.

 A closed term $\overline{\varphi}$ is interpreted as the name of the sentence $\varphi$. 

 Compound closed terms are interpreted as names of appropriate objects of the language in the 
 standard way.

 Sentences will be of  primary concern and the mechanism of refering to them will be the following. 
In the metalanguage we will use Greek letters $\varphi$, $\psi$, $\ldots$, for
 variables in sentences. So, in the language we can refer to any sentence $\varphi$ 
by a closed term  $\overline{\varphi}$. To express that something is true for all sentences,
for example that from  truth in the primary semantics follows  truth in the final semantics,
we will simply say that for every sentence $\varphi$  
 $T(\bar{\varphi})\ra \varphi$ is a true sentence of the language.  
 
 And  one more detail. Because of the uniformity of  notation, the use of the simbol 
$\bar{\ }$ is threefold. Basic use is in the construction of the term $\overline{\varphi}$
by which we refere to the sentence $\varphi$. For example, $\overline{1=1}$ 
is a name in the language $LT$ of the sentence  1=1 of $LT$. The second use is in the construction
of sentence constants by which we achieve selfreference. For example, in the expression
$\overline{L}$ it hasn't a basic use because it is not a name of the sentence $L$. Namely, there
is no such  sentence in the language $LT$. The sign $L$ we can eventually understand as
a metalanguage name for the sentence, which is, in the language named by the sentence constant
$\overline{L}$, and it is, by the previous, the sentence $T(\overline{\neg L})$. The third use 
of the simbol is in a sense of an operator which to any object $a$ of a domain of the original
language $L$ associates its name   $\overline{a}$ in the language. Contrary to the previous uses,
that which is put ``under the dash'' generally is not  an expression of the language, but an object external to
the language.

  What follows is a description of the  {\bf primary semantics}, that is  the description of the primary
 truth valuation of  sentences of the language $LT$.  Because of the
 assumption that every object $a\in DLT$ has its name $\overline{a}$ we may consider only
sentences. Valuations of arbitrary formulae and terms will be introduced later. Conditions
for the truth valuation $I_{c}$ of the sentences are  classical  together with the classical
condition for the truth predicate $T$, but what is rejected is the classical assumption that it is a
total function, defined for every sentence. Among all  functions of the kind we will select the one which is
on its domain  unique (let's remember  the possibility of multiple valuations is considered
as a failure of the classical procedure), and between all those functions we will select the 
maximal one, because we accept every success of the truth value determination. So we define the
{\bf classical truth value function} $I_{c}$ of the language  $LT$ as a partial function
$I_{c}:ST\leadsto \{\top ,\bot\}$ which obeys the following:

\begin{enumerate}
\item On atomic sentences which begin with predicates of the language $L$ values of $I_{c}$
coincide with truth values of the sentences in the language $L$ interpreted over extended domain
 $DLT$, on atomic sentences of the form  $S(\overline{a})$ it gives truth
 ($\top$) if $a$ is a sentence, otherwise falsehood  ($\bot$), and on  atomic sentences
of the form $T(\overline{a})$ where $a$ isn't a sentence it gives a falsehood. 

\item classical conditions:

      \begin{enumerate}
\item $\displaystyle I_{c}(\neg\varphi )=\left\{\begin{array}{ll}
                                             \top & for\ I_{c}(\varphi )=\bot \\
                                              \bot & for\ I_{c}(\varphi )=\top 
                                             \end{array}\right. $
\item $\displaystyle I_{c}(\varphi \wedge \psi )=\left\{\begin{array}{ll}
                                             \top & for\ I_{c}(\varphi )=\top\ and\  I_{c}(\psi )=\top\
                                               (both\ are\ true)\\
                                              \bot &  for\ I_{c}(\varphi )=\bot\ or\  I_{c}(\psi )=\bot\
                                               (at\ leat\ one\ is\ false) 
                                             \end{array}\right. $

\item $\displaystyle I_{c}(\varphi \vee \psi )=\left\{\begin{array}{ll}
                                             \top & for\ I_{c}(\varphi )=\top\ or\  I_{c}(\psi )=\top\
                                               (at\ least\ one\ is\ true)\\
                                              \bot &  for\ I_{c}(\varphi )=\bot\ and\  I_{c}(\psi )=\bot\
                                               (both\ are\ false) 
                                             \end{array}\right. $

\item $\displaystyle I_{c}(\varphi \ra \psi )=\left\{\begin{array}{ll}
                                             \top & for\ I_{c}(\varphi )=\bot\ or\  I_{c}(\psi )=\top\
                                                   \\
                                              \bot &  for\ I_{c}(\varphi )=\top\ and\  I_{c}(\psi )=\bot\
                                               
                                             \end{array}\right. $

\item $\displaystyle I_{c}(\varphi \lra \psi )=\left\{\begin{array}{ll}
                                             \top & for\ I_{c}(\varphi )=I_{c}(\psi )\
                                               (both\ are\ true\ or\ both\ are\ false)\\
                                              \bot &  for\ I_{c}(\varphi )\neq I_{c}(\psi )\
                                               (one\ is\ true\ and\ another\ is\ false) 
                                             \end{array}\right. $

\item $\displaystyle I_{c}(\forall x\varphi (x))=\left\{\begin{array}{ll}
                                             \top & if\ \forall a\in DLT\ I_{c}(\varphi (\bar{a}))=\top \\
                                              \bot & if\ \exists a\in DLT\ I_{c}(\varphi (\bar{a}))=\bot 
                                             \end{array}\right. $

\item $\displaystyle I_{c}(\exists x\varphi (x))=\left\{\begin{array}{ll}
                                             \top & if\ \exists a\in DLT\ I_{c}(\varphi (\bar{a}))=\top \\
                                              \bot & if\ \forall a\in DLT\ I_{c}(\varphi (\bar{a}))=\bot
                                             \end{array}\right. $

\end{enumerate}  

\item  classical condition on the truth predicate:

    $\displaystyle I_{c}(T(\bar{\varphi}))=\left\{\begin{array}{ll}
                                             \top & for\ I_{c}(\varphi )=\top \\
                                              \bot & for\ I_{c}(\varphi )=\bot 
                                             \end{array}\right. $
   
\item  uniqueness on the domain:

   If there is a function  $I:ST\leadsto \{\top ,\bot\}$ which obeys all three previous conditions, 
then  for every sentence $\varphi\in Dom(I_{c})\cap Dom(I_{c})$ $I(\varphi )= I_{c}(\varphi )$.

\item  maximality:

  For every function $I:ST\leadsto \{\top ,\bot\}$ which obeys all previous conditions
 
   $Dom(I)\subseteq Dom(I_{c})$.

\end{enumerate}

  From the definition  uniqueness of  such function easily follows. If there were two
such functions according to  the last condition they would have the same domain, and by the fourth condition
they would coincide on it, so they would be eqal. Later, the existence of  such function will be 
proved.

 The concept of truth values of sentences is extended to arbitrary formulae in a standard way ---
by fixing meanings of variables. The function $v:Var\longra DLT$ (where $Var$ is a set of
variables of the language) which determines meanings of variables will be called a {\bf valuation}.
In a given valuation $v$ a formula $\varphi (x_{1}, x_{2}, \ldots ,x_{n})$ with free 
variables $x_{1}, x_{2}, \ldots ,x_{n}$ is considered true  $\lra$ the associated sentence
$\varphi (\overline{v(x_{1})}, \overline{v(x_{2})}, \ldots ,\overline{v(x_{n})})$ is true,
and false $\lra$ the associated sentence is false. Also, in a given valuation $v$ 
a term $t(x_{1}, x_{2}, \ldots ,x_{n})$ with free variables $x_{1}, x_{2}, \ldots ,x_{n}$
is considered to denote the same as the closed term 
$t (\overline{v(x_{1})}, \overline{v(x_{2})}, \ldots ,\overline{v(x_{n})})$.

 To expose better the structure of the classical truth value function we will extend it to a total
function in a way that we will associate the third value $\mid$ (the undetermined) to the sentences 
on which it isn't defined. This function will be called the {\bf thruth value function of the primary 
semantics} $I_{p}:ST\longrightarrow \{\top ,\bot ,\mid ,\}$:

$\displaystyle I_{p}(\varphi )=\left\{\begin{array}{ll}
                                             I_{c}(\varphi ) & for\ \varphi\in Dom(I_{c}) \\
                                              \mid & otherwise
                                             \end{array}\right. $

 A set of sentences on which $I_{p}$ gaines  classical truth values $\top$ and $\bot$ will be called
its {\bf domain of determination}  $DDI_{p}$.

 From the definition we can easily find truth conditions (truth tables) of sentence constructions 
for the function. If arguments of the construction are
 classical  ($\top$ i $\bot$), then the value is classical too, given by the classical
conditions. If some arguments have a value the undetermined ($\mid$), then we investigate if
this failure propagates to the determination of the value of the construction on the classical 
conditions. If this is the case, then the value is also equal to the undetermined, and if it is not
the case, the value is the classical one. For example, the value of the sentence $\varphi \wedge \psi$
for $\varphi$ undetermined and $\psi$ false is false because on the classical conditions it is 
sufficient that at least one sentence is false (here it is $\psi$) for the whole sentence to be false. But if
$\psi$ is true then the truth value of the compound sentence essentialy depends on a truth value
of $\varphi$. According to the classical conditions, if $\varphi$ is true then the conjuction is also true, 
and if $\varphi$ is false than it is false, too. But $\varphi$ is undetermined, so the failure propagates
trough the conjuction which  is therefore undetermined, too. In  such way the following 
{\bf truth value conditions of the primary semantics} $I_{p}$ are given:
on

   \begin{enumerate}
\item  $\displaystyle I_{p}(\neg\varphi )=\left\{\begin{array}{ll}
                                             \top & for\ I_{p}(\varphi )=\bot \\
                                              \bot & for\ I_{p}(\varphi )=\top \\
                                              \mid & otherwise
                                             \end{array}\right. $

\vspace{1cm}

  \begin{tabular}{c|c}
 $\varphi$ & $\neg\varphi$ \\
\hline
     $\top$            & $\bot$  \\
     $\bot$            & $\top$  \\
     $\mid$            & $\mid$  
     
\end{tabular}

\vspace{1cm}

\item $\displaystyle I_{p}(\varphi \wedge \psi )=\left\{\begin{array}{ll}
                                             \top & for\ I_{p}(\varphi )=\top\ and\  I_{p}(\psi )=\top\
                                               (both\ are\ true)\\
                                              \bot &  for\ I_{p}(\varphi )=\bot\ or\  I_{p}(\psi )=\bot\
                                               (at\ least\ one\ is\ false)\\
                                              \mid & otherwise 
                                             \end{array}\right. $

\vspace{1cm}

  \begin{tabular}{c|ccc}
 $\varphi$ $\backslash$   $\psi$ & $\top$ &  $\bot$ & $\mid$ \\
\hline
                          $\top$ & $\top$ &  $\bot$ & $\mid$ \\
                          $\bot$ & $\bot$ &  $\bot$ & $\bot$ \\
                          $\mid$ & $\mid$ &  $\bot$ & $\mid$
     
\end{tabular}

\vspace{1cm}   

\item $\displaystyle I_{p}(\varphi \vee \psi )=\left\{\begin{array}{ll}
                                             \top & for\ I_{p}(\varphi )=\top\ or\  I_{p}(\psi )=\top\
                                               (at\ least\ one\ is\ true)\\
                                              \bot &  for\ I_{p}(\varphi )=\bot\ and\  I_{p}(\psi )=\bot\
                                               (both\ are\ false)\\
                                               \mid & otherwise
                                             \end{array}\right. $

\vspace{1cm}

 \begin{tabular}{c|ccc}
 $\varphi$ $\backslash$   $\psi$ & $\top$ &  $\bot$ & $\mid$ \\
\hline
                          $\top$ & $\top$ &  $\top$ & $\top$ \\
                          $\bot$ & $\top$ &  $\bot$ & $\mid$ \\
                          $\mid$ & $\top$ &  $\mid$ & $\mid$
     
\end{tabular}

\vspace{1cm}

\item $\displaystyle I_{p}(\varphi \ra \psi )=\left\{\begin{array}{ll}
                                             \top & for\ I_{p}(\varphi )=\bot\ or\  I_{p}(\psi )=\top\
                                                   \\
                                              \bot & for\ I_{p}(\varphi )=\top\ and\  I_{p}(\psi )=\bot\
                                                    \\
                                             \mid & otherwise
                                             \end{array}\right. $

\vspace{1cm}

 \begin{tabular}{c|ccc}
 $\varphi$ $\backslash$   $\psi$ & $\top$ &  $\bot$ & $\mid$ \\
\hline
                          $\top$ & $\top$ &  $\bot$ & $\mid$ \\
                          $\bot$ & $\top$ &  $\top$ & $\top$ \\
                          $\mid$ & $\top$ &  $\mid$ & $\mid$
     
\end{tabular}

\vspace{1cm}

\item $\displaystyle I_{p}(\varphi \lra \psi )=\left\{\begin{array}{ll}
                                             \top & for\ I_{p}(\varphi )=I_{p}(\psi )\neq \mid \ 
                                               (both\ are\ true\ or\ both\ are\ false)\\
                                              \bot &  for\ I_{p}(\varphi )\neq I_{p}(\psi )\ and\ none\
                                               value\ is\ \mid\\
                                              \ & (one\ is\ true\ and\ the\ other\ is\ false) \\
                                              \mid & otherwise
                                             \end{array}\right. $

\vspace{1cm}

 \begin{tabular}{c|ccc}
 $\varphi$ $\backslash$   $\psi$ & $\top$ &  $\bot$ & $\mid$ \\
\hline
                          $\top$ & $\top$ &  $\bot$ & $\mid$ \\
                          $\bot$ & $\bot$ &  $\top$ & $\mid$ \\
                          $\mid$ & $\mid$ &  $\mid$ & $\mid$
     
\end{tabular}

\vspace{1cm}

\item $\displaystyle I_{p}(\forall x\varphi (x))=\left\{\begin{array}{ll}
                                             \top & if\ \forall a\in DLT\ I_{p}(\varphi (\bar{a}))=\top \\
                                              \bot & if\ \exists a\in DLT\ I_{p}(\varphi (\bar{a}))=\bot \\
                                              \mid & otherwise
                                             \end{array}\right. $

\item $\displaystyle I_{p}(\exists x\varphi (x))=\left\{\begin{array}{ll}
                                             \top & if\ \exists a\in DLT\ I_{p}(\varphi (\bar{a}))=\top \\
                                              \bot & if\ \forall a\in DLT\ I_{p}(\varphi (\bar{a}))=\bot \\
                                             \mid  & otherwise
                                             \end{array}\right. $

\end{enumerate}

Let's note that all connectives and quantors except $\lra$ preserve their classical meaning.

 These conditions on a threevalued truth function are known in  literature as {\bf Strong
Kleene threevalued semantics} (see for example \cite{gup}). Usually, it is interpreted
as a semantics of a success of parallel algorithms or as a semantics of truth value 
investigations of sentences in a sense that sentences which hasn't {\em yet} a truth
value are declared as undetermined. Here it is interpreted as the classical procedure
 of truth value determination extended by the propagation of its own failure.

 Concerning the tuth predicate the classical determination of truth value of  $T(\bar{\varphi})$
fails precisely when determination of truth value of $\bar{\varphi}$ fails. So, $I_{p}$ obeys
the following

$$I_{p}(T(\bar{\varphi} ))=I_{p}(\varphi ) $$

 The function which obeys the additional condition is called a {\bf fixed point of Strong Kleene
semantics}.

 From the uniqueness condition on $I_{c}$ on its domain it follows an appropriate uniqueness 
condition on
 $I_{p}$ on the domain of its determination $DD(I_{p})$. Namely, for every fixed point $I$ if the
sentence $\varphi$ belongs to $DD(I)\cap DD(I_{p})$ (both valuations have a determined value on it)
then $I(\varphi )=I_{p}(\varphi )$.  It is easy to prove that it is equivalent to the following
condition on compatibility with other fixed points:

 For every fixed point  $I$ of Strong Kleene semantics  it is true that

\begin{enumerate}
\item $I_{p}(\varphi )=\top\ \ra\ I(\varphi )=\top\ or\ I(\varphi )=\mid$ 
\item $I_{p}(\varphi )=\bot\ \ra\ I(\varphi )=\bot\ or\ I(\varphi )=\mid$ 
\end{enumerate}

  Such fixed point is called  an {\bf intrinsic} point.
    
 The maximality condition of the function $I_{c}$ entails the maximality condition on  $I_{p}$
 --- for every other intrinsic fixed point $I$ $DD(I)\subseteq DD(I_{p})$.

  Therefore $I_{p}$ is a maximal intrinsic fixed point of the Strong Kleene semantics.
It is well known  (for details see for example \cite {gup}) that there is a unique such point and 
it entails that there is a unique classical truth value function $I_{c}$. Namely, for such $I_{p}$
we can define  $I_{c}:ST\leadsto \{\top ,\bot\}$ such that $D(I_{c})=DD(I_{p})$ and for every
$\varphi\in D(I_{c})$ $I_{c}(\varphi )=I_{p}(\varphi )$. It is easy to see that $I_{c}$ satisfies
all conditions on a classical truth value function. So it is proved

 {\bf Theorem}. There is a unique classical truth value function $I_{c}$ of the language $LT$.

 With this result the primary semantics of the language $LT$ is completely determined. Because 
in its threevalued formulation it is precisely the maximal intrinsic fixed point of the Strong Kleene
semantics this analysis of the classical procedure and its failures gives an argument for the choice
between various fixed point of various threevalued semantics.

  The {\bf Final semantics} $I_{f}$ of the language  $LT$ is achieved, as it has already been described 
in the first section, by taking  the primary semantics  for its subject.
The truth predicate will talk again about truth value of sentences in the primary semantics, but
now in the frame of the final semantics. Therefore, all other semantical specifications remain the
same like in the primary semantics, except for truth values  of its atomic sentences of 
the form $T(\bar{\varphi} )$ which now has an external specification:

$\displaystyle I_{f}(T(\bar{\varphi}))=\left\{\begin{array}{ll}
                                             \top & for\ I_{p}(\varphi )=\top \\
                                              \bot & otherwise
                                             \end{array}\right. $

Now the function  $I_{f}$ has given values $\top$ or $\bot$ on all atomic sentences. As it obeys
all classical conditions on truth values of compound sentences it is a classical total twovalued
truth function. By the recursion principle on the sentence structure there is a unique such function
$I_{f}:ST\longrightarrow \{\top ,\bot \}$.

From the following definitions it is clear that by the predicate $T$ we can describe the remaining
truth values of the primary semantics: 

$F(\bar{\varphi})\ \lra\ T(\overline{\neg\varphi})$

$U(\bar{\varphi})\ \lra\ \neg F(\bar{\varphi})\wedge\neg T(\bar{\varphi})$

 It is also convenient to introduce a predicate ``to have a determinate truth value''

$D(\bar{\varphi})\ \lra\  F(\bar{\varphi})\vee T(\bar{\varphi})$

 To gain  better insight in expressive power of the final semantics some sentences, which are
 true in it, will be
listed. The proofs are not given because they are straightforward.

 First of all,  in the final semantics sentences which express its consistency are true. Namely,
for every sentence $\varphi$ it is true

$\neg (T(\bar{\varphi}) \wedge F(\bar{\varphi}))$

 The following truths are  direct descriptions of truth tables of Strong Kleene semantics:

\begin{enumerate}
\item denial:
   \begin{enumerate}
   \item $T(\overline{\neg\varphi})\ \lra\ F(\bar{\varphi})$
   \item $F(\overline{\neg\varphi})\ \lra\ T(\bar{\varphi})$
   \item $U(\overline{\neg\varphi})\ \lra\ U(\bar{\varphi})$
   \end{enumerate}
\item conjuction:
   \begin{enumerate}
   \item $T(\overline{\varphi\wedge\psi})\ \lra\ T(\bar{\varphi})\wedge T(\bar{\psi})$
   \item $F(\overline{\varphi\wedge\psi})\ \lra\ F(\bar{\varphi})\vee F(\bar{\psi})$
   \item $U(\overline{\varphi\wedge\psi})\ \lra\ (T(\bar{\varphi})\wedge U(\bar{\psi}))
         \vee (U(\bar{\varphi})\wedge T(\bar{\psi}))\vee (U(\bar{\varphi})\wedge U(\bar{\psi}))$
   \end{enumerate}
\item disjunction:
   \begin{enumerate}
   \item $T(\overline{\varphi\vee\psi})\ \lra\ T(\bar{\varphi})\vee T(\bar{\psi})$
   \item $F(\overline{\varphi\vee\psi})\ \lra\ F(\bar{\varphi})\wedge F(\bar{\psi})$
   \item $U(\overline{\varphi\vee\psi})\ \lra\ (F(\bar{\varphi})\wedge U(\bar{\psi}))
         \vee (U(\bar{\varphi})\wedge F(\bar{\psi}))\vee (U(\bar{\varphi})\wedge U(\bar{\psi}))$
   \end{enumerate}
\item conditional:
   \begin{enumerate}
   \item $T(\overline{\varphi\ra\psi})\ \lra\ F(\bar{\varphi})\vee T(\bar{\psi})$
   \item $F(\overline{\varphi\ra\psi})\ \lra\ T(\bar{\varphi})\wedge F(\bar{\psi})$
   \item $U(\overline{\varphi\ra\psi})\ \lra\ (T(\bar{\varphi})\wedge U(\bar{\psi}))
         \vee (U(\bar{\varphi})\wedge F(\bar{\psi}))\vee (U(\bar{\varphi})\wedge U(\bar{\psi}))$
   \end{enumerate}
\item biconditional:
     \begin{enumerate}
   \item $T(\overline{\varphi\lra\psi})\ \lra\ (T(\bar{\varphi})\wedge T(\bar{\psi}))
                                               \vee (F(\bar{\varphi})\wedge F(\bar{\psi}))$
   \item $F(\overline{\varphi\lra\psi})\ \lra\ (T(\bar{\varphi})\wedge F(\bar{\psi}))
                                               \vee (F(\bar{\varphi})\wedge T(\bar{\psi}))$
   \item $U(\overline{\varphi\lra\psi})\ \lra\ U(\bar{\varphi})\vee U(\bar{\psi})$
         
   \end{enumerate}
\item universal quantification: 
    \begin{enumerate}
   \item  $T(\overline{\forall x\varphi (x)})\lra 
          \forall x T(\overline{\varphi (x)})$
   \item $F(\overline{\forall x\varphi (x)})\lra \exists x F(\overline{\varphi (x)})$
   \item $U(\overline{\forall x\varphi (x)})\lra \neg\exists x F(\overline{\varphi (x)})
          \wedge\exists x U(\overline{\varphi (x)})$
   \end{enumerate}
\item existential quantification: 
    \begin{enumerate}
   \item  $T(\overline{\exists x\varphi (x)})\lra \exists x T(\overline{\varphi (x)})$
   \item $F(\overline{\exists x\varphi (x)})\lra \forall x F(\overline{\varphi (x)})$
   \item $U(\overline{\exists x\varphi (x)}) \lra \neg\exists x T(\overline{\varphi (x)})\wedge
    \exists x U(\overline{\varphi (x)})$
   \end{enumerate} 
\end{enumerate}

  The iteration of the truth predicate is not interesting because the following is true:

\begin{enumerate}
   \item $T(\overline{T(\overline{\varphi})})\ \lra\ T(\bar{\varphi})$
   \item $F(\overline{T(\overline{\varphi})})\ \lra\ F(\bar{\varphi})$
   \item $T(\overline{F(\overline{\varphi})})\ \lra\ F(\bar{\varphi})$
   \item $F(\overline{F(\overline{\varphi})})\ \lra\ T(\bar{\varphi})$
   \item $U(\overline{T(\overline{\varphi})})\ \lra\ U(\bar{\varphi})$ etc.
   \end{enumerate}

  Previous rules reduce the investigation of the truth value of a sentence $\varphi$
in the primary semantics, that is to say an investigation of  truth values of $T(\bar{\varphi})$
and $F(\bar{\varphi})$ in the final semantics, to an investigation of truth values in the final
semantics of atomic sentences of the form $T(t)$ and $F(t)$, where $t$ is a name of 
an object which isn't a sentence, or $t=\bar{\psi}$ where  $\psi$  is an atomic sentence which
doesn't begin with predicate $T$ or $F$, or it is a sentence constant. In first two cases
the answer is simple:

 If $t$ isn't a name of a sentence then 

$\neg T(t)\wedge \neg F(t)$

 If  $\psi$ is an atomic sentence which doesn't began with $T$ or $F$ then 

\begin{enumerate}
\item $T(\bar{\psi})\lra\psi$
\item $F(\bar{\psi})\lra\neg\psi$
\end{enumerate}

  Therefore, what remains is to determine truth values of sentences $T(C)$ and $F(C)$ where  $C$ is a 
sentence constant. It is the most interesting part of the language  because 
 selfreferent sentences are constructed by the
sentence constants . Principally we can determine
truth values of such sentences by the analysis of the primary semantics, but it is interesting 
to see in what amount  selfreference and intuitive argumentation leading to paradoxes in the
classical language can be reproduced in the language $LT$. In an intuitive argumentation a transition
from  assumption about truth value of a sentence to acception or rejection of what it says
is a crucial step. At  first sight we can desribe it in the final semantics using the sentences

$T(C)\ra C$ and $F(C)\ra \neg C$

  But there is a technical problem that $C$ is a sentence constant which names
a sentence, let's say $d(C)$, and not the sentence itself. So a correct description is

\begin{enumerate}
\item $T(C)\ra d(C)$
\item $F(C)\ra \neg d(C)$
\end{enumerate}

 For example, the sentence constant $\overline{LL}$ names 
$d(\overline{LL})=\neg T(\overline{LL})$ (the Strong  Liar) so its description in  $LT$ is

\begin{enumerate}
\item $T(\overline{LL})\ra \neg T(\overline{LL})$
\item $F(\overline{LL})\ra \neg \neg T(\overline{LL})$
\end{enumerate} 

 Using the description we can translate the intuitive argumentation in the language $LT$.
Let $T(\overline{LL})$ be true. Then, by the first description, $\neg T(\overline{LL})$ is
true and it is a contradiction. If we assume $F(\overline{LL})$ then, by the second description,
$\neg \neg T(\overline{LL})$, that is  $T(\overline{LL})$ is true and it is also in a contradiction
with the statement of consistency of the primary semantics $\neg (T(\bar{\varphi}) \wedge F(\bar{\varphi}))$.
So, $U(\overline{LL})$. Particulary, it means that $\neg T(\overline{LL})$. Therefore, we showed in the
language $LT$ that the Strong Liar is undetermined in the primary semantics and true in the final semantics.

In the same way, every paradox  which leads to a contradiction in a classical semantics can be translated
in an argumentation in the final semantics which states truth values of a sentence in the primary
and final semantics. But, as it has already been shown in the previous section in an informal way, such
description is not sufficient to state truth values of selfereference sentences, which don't lead to
a contradiction, but permit one truth valuation (as the Logician) or more (as the Truthteller). For
example, for the Truthteller the description is

 \begin{enumerate}
\item $T(\overline{I})\ra  T(\overline{I})$
\item $F(\overline{I})\ra \neg T(\overline{I})$
\end{enumerate}

 But it is true for every sentence and we can deduce nothing about the Truthteller. 
In such cases it is necessary to look at the semantical graph and state the primary valuation 
of the sentence
$\varphi$ to know in the language $LT$ what is a truth value of $T(\overline{\varphi})$ and of $\varphi$.

We will display  some other principles which talk about the truth predicate  $T$ (and $F$). Of 
course, we know that Tarski's schema  $T(\overline{\varphi})\lra \varphi$ for every $\varphi$
is not valid (\cite{tar}). Here, it is a consequence of the fact that $T$ is not a truth predicate
 for the final, but for the primary semantics of the language. But it is true that the final semantics
is an extension of the primary one by the description of its failures. Everything true in the primary semantics
is true in the final semantics, and everything false in the primary one is false in the final one, that is 
for every sentence $\varphi$ it is true:

\begin{enumerate}
\item $T(\overline{\varphi})\ra \varphi$
\item $F(\overline{\varphi})\ra \neg\varphi$
\end{enumerate}

 Of course, for sentences which have a definite truth value in the primary semantics the converse
is also true:

 $D(\overline{\varphi})\ \ra\ (T(\overline{\varphi})\lra \varphi) \wedge (F(\overline{\varphi})\lra \neg\varphi)$

 For $\varphi _{1}$ and $\varphi _{2}$ logically equivalent sentences in classical logic it is true
that they are logically equivalent in the primary semantics, so it is true in the final semantics:

\begin{enumerate}
   \item $T(\overline{\varphi _{1}})\ \lra\ T(\overline{\varphi _{2}})$
   \item $F(\overline{\varphi _{1}})\ \lra\ F(\overline{\varphi _{2}})$
   \item $U(\overline{\varphi _{1}})\ \lra\ U(\overline{\varphi _{2}})$
   \end{enumerate}

\bibliography{sinth}
 \addcontentsline{toc}{section}{Literature}
{\em   
Boris \v{C}ulina \\
Department of Mathematics\\
Faculty of Mechanical Engineering and Naval Architecture\\
University of Zagreb\\
Ivana Lu\kj i\tj a 5\\
10 000 Zagreb\\
CROATIA}
 \end{document}